\documentclass{amsart}

\usepackage{latexsym}
\usepackage{amsmath,amssymb,stmaryrd}
\usepackage{amscd,amsxtra}
\usepackage{amsthm}

\theoremstyle{plain} 
\newtheorem{thm}{\bf Theorem}[section]
\newtheorem{pro}[thm]{\bf Proposition}
\newtheorem{lem}[thm]{\bf Lemma}

\DeclareMathOperator{\ann}{Ann}

\title{On the Frobenius functor and colon ideals}
\author{Wenliang Zhang}
\address{Department of Mathematics, University of Minnesota, Minneapolis, MN 55455}
\email{wlzhang@math.umn.edu}

\begin{document}

\begin{abstract}
In this paper we study the commutativity of the Frobenius functor and the colon operation of two ideals for Noetherian rings of positive characteristic $p$. New characterizations of regular rings and local UFDs are given.
\end{abstract}
\maketitle

\section{Introduction}
In this short note, we study the commutativity of the Frobenius functor and the colon operation of two ideals for Noetherian rings of positive characteristic $p$. To this end, we introduce the following condition
\newtheorem*{ci}{Condition ($C_i$)}
\label{ci}
\begin{ci}
$((x_1,...,x_i):y)^{[q]}=((x^q_1,...x^q_i):y^q)$ for all elements $x_1,...,x_i,y\in R$ and all $q=p^e$, where $R$ is a Noetherian ring of positive characteristic $p$ (for $i=0$, we set $(x_1,\dots,x_i)=0$).
\end{ci}
We say that $R$ satisfies $C_{\infty}$ if $R$ satisfies $C_i$ for all $i\geq 0$.\par
A natural question is, can one characterize rings satisfying the condition ($C_i$)?\par
It is easy to see that a Noetherian local ring of positive characteristic $p$ satisfies condition ($C_0$) if and only if it is a domain. In Section 2, we show that the condition ($C_{\infty}$) characterizes regular rings of positive characteristic. And in Section 3, we prove that the condition ($C_1$) characterizes local UFDs.
\section{A characterization of regular rings}
Characterizing regular rings has a long and rich history. For example, a Noetherian local ring is regular if and only if it has finite global dimension (due to Serre, Theorem 19.2 in \cite{m}). In particular, if the characteristic is positive, there is Kunz's characterization of regular local rings.  

\begin{thm}{\bf (Theorem 2.1, Theorem 3.3 in \cite{k})}
\label{kunz}
Let $(R,\mathfrak{m})$ be a Noetherian local ring of positive characteristic $p>0$. Then the following are equivalent
\begin{enumerate}
\item $R$ is regular 
\item The Frobenius functor is exact
\item For some $q=p^e$ with $e\geq 1$ (equivalently for all $q$), $\lambda_R(R/\mathfrak{m}^{[q]})=q^d$.
\end{enumerate} 
\end{thm}

In this section, we will prove the following characterization of regular rings. 

\begin{thm}
\label{regular-local}
Let $(R,\mathfrak{m})$ be a Noetherian local ring of positive characteristic $p$. Then the following are equivalent 
\begin{enumerate}
\item $R$ is a regular ring.
\item $(I^{[q]}:J^{[q]})=(I:J)^{[q]}$ for any pair of ideals $I$ and $J$ and all $q=p^e$, $e\geq 0$.
\item $(\mathfrak{a}^{[q]}:x^q)=(\mathfrak{a}:x)^{[q]}$ for all $\mathfrak{m}$-primary ideals $\mathfrak{a}$, all elements $x\in R$, and all $q=p^e$ for $e\geq 0$.
\end{enumerate}
\end{thm}

A lemma is needed before one can prove this characterization.

\begin{lem}
\label{lambda}
Let $(R,\mathfrak{m})$ be a Noetherian local ring of characteristic $p>0$. Assume that $(\mathfrak{a}^{[q]}:x^q)=(\mathfrak{a}:x)^{[q]}$ for all $\mathfrak{m}$-primary ideals $\mathfrak{a}$, all element $x\in R$, and all $q=p^e$ for $e\geq 1$. Let $I\subseteq J$ be 2 ideals with $I$ $\mathfrak{m}$-primary. Then 
$$\lambda(R/I^{[q]})=\lambda(J/I)\cdot\lambda(R/\mathfrak{m}^{[q]})+\lambda(R/J^{[q]}).$$
\end{lem}
\begin{proof}[Proof]
Set $l=\lambda(J/I)$. Choose any filtration of $I\subseteq J\subseteq R$,
$$I=J_0\subsetneq J_1\subsetneq J_2\subsetneq\cdots\subsetneq J_l=J\subseteq R,$$
such that $J_i/J_{i-1}\cong R/\mathfrak{m}$, for all $i=1,...,l$. That is to
say, $J_i=(J_{i-1},x_i)$ for some $x_i\in J_i$ such that $(J_{i-1}:x_i)=\mathfrak{m}$. (One can set $x_i$ to be a preimage in $J_i$ of any generator of $J_i/J_{i-1}$. $(J_{i-1}:x_i)=\mathfrak{m}$ follows from the short exact sequence $0\to R/(J_{i-1}:x_i)\to R/J_{i-1}\to R/J_i\to 0$.) For all 
$q$, there is a corresponding filtration of $I^{[q]}\subseteq J^{[q]}\subseteq R$,
\begin{equation}
\label{filtration}
I^{[q]}=J^{[q]}_0\subseteq J^{[q]}_1\subseteq \cdots\subseteq J^{[q]}_l=J^{[q]}\subseteq R,
\end{equation}
where $J^{[q]}_i/J^{[q]}_{i-1}\cong R/(J^{[q]}_{i-1}:x^q_i)$ (it follows from the short exact sequence $0\to R/(J^{[q]}_{i-1}:x^q_i)\to R/J^{[q]}_{i-1}\to R/J^{[q]}_i\to 0$). But
$(J^{[q]}_{i_1}:x^q_i)=(J_{i-1}:x_i)^{[q]}=\mathfrak{m}^{[q]}$, for each $i=1,...,l$, hence $\lambda(R/I^{[q]})=\lambda(J/I)\cdot\lambda(R/\mathfrak{m}^{[q]})+\lambda(R/J^{[q]})$ (from the additivity of $\lambda$ and the filtration (\ref{filtration})).
\end{proof}

\textbf{Remark}. In Lemma \ref{lambda}, if one does not have the assumption ``$(\mathfrak{a}^{[q]}:x^q)=(\mathfrak{a}:x)^{[q]}$ for all $\mathfrak{m}$-primary ideals $\mathfrak{a}$, all element $x\in R$, and all $q=p^e$ for $e\geq 0$", then the inequality 
$$\lambda(R/I^{[q]})\leq\lambda(J/I)\cdot\lambda(R/\mathfrak{m}^{[q]})+\lambda(R/J^{[q]})$$
always holds (cf. Lemma 2.1 in \cite{hy} or Proposition 5.2.1 in \cite{h}).

\begin{proof}[Proof of Theorem \ref{regular-local}]
(1)$\Rightarrow$(2) is a well-known consequence from Theorem 1.1(2); (2)$\Rightarrow$(3) is trivial. It remains to prove that (3)$\Rightarrow$(1).\par
Set $J=R$ and $I=\mathfrak{m}^{[p]}$ in Lemma 1.3. We have 
$$\lambda(R/\mathfrak{m}^{[pq]})=\lambda(R/\mathfrak{m}^{[p]})\cdot\lambda(R/\mathfrak{m}^{[q]}).$$
Dividing by $q^d$, $d=\dim(R)$, and taking limit over $q\to \infty$, we have
$$p^d\cdot e_{HK}(R)=\lambda(R/\mathfrak{m}^{[p]})\cdot e_{HK}(R),$$
where $e_{HK}(R)$ is the Hilbert-Kunz multiplicity of $R$ and its existence is guaranteed by \cite{m1}, hence, $\lambda(R/\mathfrak{m}^{[p]})=p^d$, Theorem 1.1 completes the proof.
\end{proof}

{\bf Remark}. It is easy to see that all three conditions (1), (2) and (3) in Theorem 1.2 `localize'. One can prove it as follows.\par 
First, we prove $((I:J)^{[q]})_W=(I_W:J_W)^{[q]}$. ``$\supseteq$", let $(\frac{x}{s})^q\in (I_W:J_W)^{[q]}$, $xJ\subseteq I,s\in W$. Then obviously, $\frac{x^q}{s^q}\in ((I:J)^{[q]})_W$. ``$\subseteq$", let $\frac{x^q}{s}\in ((I:J)^{[q]})_W$, where $x\in (I:J)$. $\frac{x^q}{s}=s^{q-1}(\frac{x}{s})^q$. But $\frac{x}{s}\in (I_W:J_W)$, $(\frac{x}{s})^q\in (I_W:J_W)^{[q]}\Rightarrow \frac{x^q}{s}=s^{q-1}(\frac{x}{s})^q\in (I_W:J_W)^{[q]}$.\par
Next, we prove $(I:J)^{[q]}=(I^{[q]}:J^{[q]})\Rightarrow (I_W:J_W)^{[q]}=(I_W)^{[q]}:(J_W)^{[q]}$. ``$\subseteq$" is obvious. ``$\supseteq$", $\frac{x}{s}\in ((I_W)^{[q]}:(J_W)^{[q]})$ iff 
$\frac{x}{s}(J_W)^{[q]}\subseteq (I_W)^{[q]}$. Then $\frac{x}{s}((\frac{j_1}{s_1})^q,...,(\frac{j_r}{s_r})^q)\subseteq (I_W)^{[q]}$ (we may assume that $(j_1,...,j_r)=J$), clearing the denominators, we have $\tilde{s}xJ^{[q]}\subseteq I^{[q]}$, then $x\in (I^{[q]}:J^{[q]})_W=((I:J)^{[q]})_W=(I_W:J_W)^{[q]}$.\par
Hence we have

\begin{thm}
\label{regular}
Let $R$ be a Noetherian ring of positive characteristic $p$. Then the following are equivalent
\begin{enumerate}
\item $R$ is a regular ring.
\item $(I^{[q]}:J^{[q]})=(I:J)^{[q]}$ for all ideals $I$ and $J$ and for all $q=p^e$, $e\geq 0$.
\item  $(\mathfrak{a}^{[q]}:x^q)=(\mathfrak{a}:x)^{[q]}$ for all ideals $\mathfrak{a}$ which are primary to some maximal ideal, all elements $x\in R$, and all $q=p^e$, $e\geq 0$.
\end{enumerate}
\end{thm}

From Theorem \ref{regular}, one can see that the condition ($C_{\infty}$) characterizes regular rings of positive characteristic.
\section{On the condition ($C_1$)}
In this section, we will show that the condition ($C_1$) characterizes local UFDs.
\begin{lem}
\label{principal-tight}
Let $R$ be a Noetherian local ring of positive characteristic $p$ and  $I$ an proper ideal of $R$. If $(I^{[q]}:r^q)=(I:r)^{[q]}$ for all $r\in R$ and $q=p^e,e\geq 1$. Then $I=I^*$.
\end{lem}
\begin{proof}[Proof]
All we need to prove is that $I^*\subseteq I$. Assume $x\in I^*$, then there exists $c\in R^o$ such that $cx^q\in I^{[q]}$ for $q\gg 0$. Hence, $c\in (I^{[q]}:x^q)=(I:x)^{[q]}$ for $q\gg 0$, which implies $c\in \bigcap_{q\gg 0}(I:x)^{[q]}$. Since $c\neq 0$, $(I:x)$ has to be $R$, i.e., $x\in I$.  
\end{proof}
Before we state and prove the key lemma (Lemma \ref{principal}), let us recall some theorems needed in the proof of Lemma \ref{principal}.

\begin{thm}{\rm (Northcott-Rees Theorem)}(cf. \cite{nr})
\label{nr}
Let $R$ be a Noetherian local ring and $I$ be an ideal of $R$. Define a function $f$ assigning each positive integer $n$ to the minimal number of generators of $I^n$. Then there exists a polynomial $H$ such that $H(n)=f(n)$ for $n\gg 0$. Moreover, $\deg(H)=l(I)+1$, where $l(I)$ is the analytic spread of $I$.
\end{thm}

\begin{pro}{\rm (Proposition 8.3.8 in \cite{hs})}
\label{integral}
Let $R$ be a Noetherian local ring and $I$ be an ideal of $R$. Then there exists an integer $n$ such that $I^n$ is integral over an ideal generated by $l(I)$ elements.
\end{pro}

\begin{thm}{\rm (Theorem 5.4 in \cite{hh1})}
\label{Briancon-Skoda}
Let $R$ be a Noetherian ring of positive characteristic $p$, and let $I$ be an ideal of positive height generated by $n$ elements. Then for every $m\in\mathbb{N}$, $\overline{I^{n+m}}\subseteq (I^{m+1})^*$. In particular, $\overline{I^n}\subseteq I^*$.
\end{thm}

\begin{lem}
\label{principal}
Let $R$ be a Noetherian local ring of positive characteristic $p$. If $((x^q):y^q)=(x:y)^{[q]}$ for all $x,y\in R$, all $q=p^e$, $e\geq 0$, then $((x):y)$ is a principal ideal for all $x,y\in R$.
\end{lem}
\begin{proof}[Proof]
First, we will prove that $R$ is a domain. Assume that $xy=0$ and $x\neq 0$. Then $yx^q=0$ for all $q$, i.e., $y\in ((0):x^q)$. Hence, $y\in \bigcap_q((0):x^q)=\bigcap_q((0):x)^{[q]}=(0)$ (since $((0):x)$ is a proper ideal), i.e., $y=0$. So, $R$ is a domain.\par 
Assume that $((x):y)$ is generated by $g_1,...,g_t$. Let $n_1,...,n_t$ be nonnegative integers such that $n_1+\cdots+n_t=q$, then $g^{n_1}_1\cdots g^{n_t}_ty^q=g^{n_1}_1y^{n_1}\cdots g^{n_t}_ty^{n_t}\in (x^q)$. Hence, $g^{n_1}_1 \cdots g^{n_t}_t\in ((x^q):y^q)=((x):y)^{[q]}$. But $((x):y)^q$ is
generated by $\{g^{n_1}_1\cdots g^{n_t}_t|\sum^t_{i=1} n_i=q\}$, $((x):y)^q=((x):y)^{[q]}$. It is obvious that the minimal nubers of generators of $((x):y)^{[q]}$ is bounded by $t$ (a constant), hence the polynomial associated to $((x):y)$ (cf. Theorem \ref{nr}) has to be a constant, i.e., has degree 0. Therefore, the analytic spread of $((x):y)$ is 1. Hence, by Proposition \ref{integral}, there exists an integer $n$ such that $((x):y)^n$ is integral over a principal ideal, say, $(a)$, where $a\in ((x):y)^n$ is a  nonzero element. Since $R$ is a domain, $(a)$ has positive height. Applying Theorem \ref{Briancon-Skoda}, one has $\overline{(a)}=(a)^*$ ($I^*\subseteq \bar{I}$ is always true all ideals $I$, Theorem 5.2 in \cite{hh1}). But Lemma \ref{principal-tight} implies that $(a)=(a)^*$ So, $(a)$ is integrally closed and hence $((x):y)^n=(a)$. $((x):y)^n$ is integral over a principal ideal, hence evrey power of $((x):y)^n$ is integral over a 
principal ideal. By the similar argument as above, each power of $((x):y)^n$ is 
principal. Therefore, the polynomial associated to $((x):y)$ is $1$. For $m\gg 0$, $((x):y)^m$ is principal, in particaul, $((x):y)^q$ is principal for $q\gg 0$. Assume $((x):y)^q=(g)$. Then there exist $r_i, s_i\in R$ such that $g^q_i=r_ig$ and $g=\sum_is_ig^q_i$. Then $g=\sum_ir_is_ig$. Since $R$ is a domain, $\sum_ir_is_i=1$. There is at least one $r_i$, say, $r_1$, which is a unit. We may assume that $((x):y)^q=(g^q_1)$. Then for any $g_i$, $g^{q-1}_1g_i\in ((x):y)^q$, i.e., there exists $t_i\in R$ such that $g^{q-1}_1g_i=t_ig_1^q$. But then $g_i=t_ig_1$, i.e., $((x):y)=(g_1)$.
\end{proof}

\begin{thm}
\label{ufd}
Let $R$ be a Noetherian commutative local ring of positive characteristic $p$. Then the following are equivalent
\begin{enumerate}
\item $R$ is a UFD.
\item $((x^q):y^q)=((x):y)^{[q]}$ for all elements $x$ and $y$ in $R$ and all $q=p^e$, $e\geq 0$.
\end{enumerate}
\end{thm}
\begin{proof}[Proof]
(1)$\Rightarrow$(2). If one of $x$ and $y$ is 0 or a unit, the conclusion is clear. we may assume that $x,y\neq 0$ and that $x$ and $y$ are not units. We can write $x=u\prod^m_{i=1} x_i$,  where $u$ is a unit and $x_i$ are prime elements. $y$ can be written as $y=v\prod^n_{j=1}y_j$, where $v$ is a unit and $y_j$ are prime elements. We will use induction on $n$ to prove that $((x^q):y^q)=((x):y)^{[q]}$.\par
$n=1$. Assume $ay^q=bx^q$, i.e., $av^qy^q_1=bu^q\prod^m_{i=1}x^q_i$. If $y_1$ divides some $x_i$, say, $x_1$, then $y_1=x_1$ (both are prime elements). Hence $a=b(\frac{u}{v}\prod^m_{i=2}x_i)^q\in ((x):y)^{[q]}$ ($y_1=x_1$ implies that $\frac{u}{v}\prod^m_{i=2}x_i\in ((x):y)$). If $y_1$ does 
not divide any $x_i$, then $y^q_1$ must divide $b$. Write $b=b'y^q_1$, then we have $a=b'(\frac{u}{v}\prod^m_{i=1}x_i)^q\in ((x):y)^{[q]}$. So, $((x^q):y^q)\subseteq ((x):y)^{[q]}$. Since the other direction is always true, $((x^q):y^q)=((x):y)^{[q]}$.\par
$n\geq 2$. Assume $ax^q=by^q$, i.e., $av^q\prod^n_{j=1}x^q_j=bu^q\prod^m_{i=1}y^q_i$. If $y_n$ divides some $x_i$, say, $x_1$, then $y_n=x_1$ and $ay^q\in (x)^{[q]}$ iff $av^q\prod^{n-1}_{j=1}y^q_j\in (x)^{[q]}$. If $y_n$ does not divide any $x_i$, then $y^q_n$ must divide $b$, again we have $ay^q\in (x)^{[q]}$ iff $av^q\prod^{n-1}_{j=1}r^q_j\in (x)^{[q]}$. Now the inductive hypotheses 
complete the proof. \par  
(2)$\Rightarrow$(1). We will prove that $R$ is a domain and that every prime ideal of height 1 is
principal.\par
From the proof of Lemma \ref{principal}, one can see that $R$ is a domain. Let $\mathfrak{p}$ be a prime ideal of height 1. Let $a$ be a nonzero element in $\mathfrak{p}$. Then $\mathfrak{p}$ is an associated prime of $(a)$. There exists $\bar{b}\in R/(a)$ such that $\ann(\bar{b})=\mathfrak{p}$. Let $b$ be any preimage of $\bar{b}$ in $R$, then $\mathfrak{p}=((a):b)$. Applying Lemma \ref{principal}, one has $\mathfrak{p}$ is principal.
\end{proof}

\textbf{Remark}. In general, let $R$ be a Noetherian integral domain, one can prove that $R$ is a UFD if and only if $((x):y)$ is pricipal for all elements $x,y\in R$.\par

\section{Acknowledgements}
The author started working on this note during his visit to University of Kansas. He would like to thank Department of Mathematics, University of Kansas for the wonderful hospitality. He also wants to thank Prof. Craig Huneke for useful discussions.

\end{document}